\documentclass{kms-j}

\theoremstyle{plain}
\newtheorem{theorem}{Theorem}[section]
\newtheorem{proposition}[theorem]{Proposition}

\theoremstyle{definition}
\newtheorem*{definition}{Definition}

\theoremstyle{remark}

\begin{document}
\title{Partition number identities\\which are true for all set of parts}

\author{Kim, Bongju}
\address{Department of Mathematics \\Pusan National University\\ Korea}
\email{kim10057@pusan.ac.kr}

\footnote{The theorems in this paper stated and proved in 2012 winter. 2010 \textit{Mathematics Subject Classification.}  11P81, 11P84, 05A17}

\begin{abstract}
Let $B$ be an infinite subset of $\mathbf{N}$. When we consider partitions of natural numbers into elements of $B$, a partition number without  a restriction of the number of equal parts can be expressed by partition numbers with a restriction $\alpha$ of the number of equal parts. Although there are many way of the expression, we prove that there exists a expression form  such that this expression form is true for all possible set $B$.  This identities comes from the partition numbers of natural numbers into $\{1,\alpha,\alpha^2,\alpha^3,\cdots\}$. Furthermore, we prove that there exist inverse forms of the expression forms. And we prove other similar identities. The proofs in this paper are constructive.
\end{abstract}

\maketitle

Let $p(n)$ be the number of partitions of $n$ into natural numbers, and $d(n)$ be the number of partitions of $n$ into distinct natural numbers. If one calculates $p(5)$ and from $d(1)$ to $d(5)$ except $d(4)$, \[p(5)=7, d(1)=1, d(2)=1, d(3)=2, d(5)=3.\]
However, $p(5)$ can de expressed by the following product and sum of $d(1)$, $d(2)$, $d(3)$, $d(5)$.
\begin{align*}
p(5)=7 & =3+1\times1+1\times1+2\times1\\
 &=d(5)+d(1)d(1)+d(2)d(1)+d(3)d(1).
\end{align*}
Next, let $p^{\mathbf{P}}(n)$ be the number of partitions of $n$ into primes, and $p^{\mathbf{P}}_1(n)$ be the number of partitions of $n$ into distinct primes. Then \[p^{\mathbf{P}}(5)=2, p^{\mathbf{P}}_1(1)=0, p^{\mathbf{P}}_1(2)=1, p^{\mathbf{P}}_1(3)=1,p^{\mathbf{P}}_1(5)=2\]
and $p^{\mathbf{P}}(5)$ can de expressed by the following product and sum of $p^{\mathbf{P}}_1(1)$, $p^{\mathbf{P}}_1(2)$, $p^{\mathbf{P}}_1(3)$, $p^{\mathbf{P}}_1(5)$.
\begin{align*}
p^{\mathbf{P}}(5)=2& =2+0\times0+1\times0+1\times0\\ &=p^{\mathbf{P}}_1(5)+p^{\mathbf{P}}_1(1)p^{\mathbf{P}}_1(1)+p^{\mathbf{P}}_1(2)p^{\mathbf{P}}_1(1)+p^{\mathbf{P}}_1(3)p^{\mathbf{P}}_1(1).
\end{align*}
Finally, let $P^{O}(n)$ be the number of partitions of $n$ into odd numbers, and $p^{O}_1(n)$ be the number of partitions of $n$ into distinct odd numbers. Then \[p^{O}(5)=3, p^{O}_1(1)=1, p^{O}_1(2)=0, p^{O}_1(3)=1, p^{O}_1(5)=1\]
and
\begin{align*}p^{O}(5)=3&=1+1\times1+0\times1+1\times1\\ &=p^{O}_1(5)+p^{O}_1(1)p^{O}_1(1)+p^{O}_1(2)p^{O}_1(1)+p^{O}_1(3)p^{O}_1(1).
\end{align*}
If one compares above three expressions, one can see that these three expression forms are same though the set of parts was changed. The numbers in the previous identity are not random but come from the all possible binary expressions of $5$ i.e. \begin{align*}
5&=5\times2^0\\
 &=1\times2^0+1\times2^2\\
 &=1\times2^0+2\times2^1\\
 &=3\times2^0+1\times2^1.
\end{align*}

\begin{definition}
Let $\mathbf{N}$ be the set of natural numbers and $\psi$ be an one to one function such that $\psi:\mathbf{N}\rightarrow \mathbf{N}$. We denote the set $\{\psi(n)\in\mathbf{N}\}$ by $A$.
$p^{A}_{\alpha}(n)$ is the number of partitions of $n$ into elements of $A$ such that the number of equal parts is less than or equals to $\alpha\in \mathbf{N}\setminus\{0\}$.
$p^{A}_{}(n)$ is the number of partitions of $n$ into elements of $\psi(\mathbf{N})$ without a restriction of the number of equal parts.
We define $p^{A}_{\alpha}(0)=p^{A}_{}(0):=1$ for all $\alpha\in\mathbf{N}\setminus\{0\}$.
\end{definition}

\section{Identities between $p^{A}_{}$ and $p^{A}_{\alpha}$}

Before considering the identities between $p^{A}_{}$ and $p^{A}_{\alpha}$, we prove the identities between $p^{A}_{}$ and $p^{A}_{1}$.

\begin{definition}
Consider all non-negative integer solutions of the indeterminate equation $n=N_{0}+2N_{1}+4N_{2}+\cdots=\sum_{i\geq 0}2^{i}N_{i}$ and denote by $(a^{n}_{11},a^{n}_{12},a^{n}_{13}\cdots)$, $(a^{n}_{21},a^{n}_{22},a^{n}_{23}\cdots)$, $\cdots$. In other word,
\begin{align*}
n&=a^{n}_{11}+2a^{n}_{12}+4a^{n}_{13}+\cdots\\
 &=a^{n}_{21}+2a^{n}_{22}+4a^{n}_{23}+\cdots\\
 &\cdots
\end{align*} Then the solution matrix of this equation is  \[A_n:=(a^{n}_{ij}).\]
\end{definition}

\begin{proposition}
Let $A_n=(a^{n}_{ij})$ be the solution matrix of $n=\sum_{i\geq 0}2^{i}N_{i}$ and let $\psi$ be an one to one function such that $\psi:\mathbf{N}\rightarrow \mathbf{N}$. Then
\[
p^{A}_{}(n)=\sum_{i\geq 1}\prod_{j\geq 1}p^{A}_{1}(a^{n}_{ij})
\]
for all $n\in\mathbf{N}$.
\end{proposition}

\begin{proof}
It is well known fact that for $|q|<1$,
\[\sum_{n\geq 0}p^{A}_{}(n)q^{n}=\prod_{n\geq 1}\frac{1}{1-q^{\psi(n)}}\]
and
\[\sum_{n\geq 0}p^{A}_{1}(n)q^{n}=\prod_{n\geq 1}(1+q^{\psi(n)})\](see \cite{1} or \cite{2}). On the other hand,
\begin{align*}
\prod_{n\geq1}(1-q^{2\psi(n)})&=\prod_{n\geq1}(1-q^{\psi(n)})\prod_{n\geq1}(1+q^{\psi(n)})\\
\prod_{n\geq1}(1-q^{4\psi(n)})&=\prod_{n\geq1}(1-q^{\psi(n)})\prod_{n\geq1}(1+q^{\psi(n)})\prod_{n\geq1}(1+q^{2\psi(n)}) \\
\vdots\\
\prod_{n\geq1}(1-q^{2^{I+1}\psi(n)})&=\prod_{n\geq1}(1-q^{\psi(n)})\prod_{i=0}^{I}\prod_{n\geq1}(1+q^{2^{i}\psi(n)})
\end{align*}
for all $I \in \mathbf{N}$. So,
\begin{align*}
1&=\lim_{I\rightarrow\infty}\prod_{n\geq1}(1-q^{2^{I+1}\psi(n)})\\
&=\prod_{n\geq1}(1-q^{\psi(n)})\prod_{i=0}^{\infty}\prod_{n\geq1}(1+q^{2^{i}\psi(n)}).
\end{align*}
Therefore,
\begin{align*}
\prod_{n\geq1}\frac{1}{1-q^{\psi(n)}}=\prod_{i\geq0}\prod_{n\geq1}(1+q^{2^{i}\psi(n)})
\end{align*}
and
\begin{align*}
\sum_{n\geq 0}p^{A}_{}(n)q^{n}
=\prod_{i\geq0}(\sum_{n\geq0}p^{A}_{1}(n)q^{2^{i}n}).
\end{align*}
If one expands the above infinite product to an infinite series, one can express $p^{A}_{}(n)$ by $p^{A}_{1}(n)$, $p^{A}_{1}(n-1)$, $\cdots$, $p^{A}_{1}(1)$ and one can see that the coefficient of $q^n$ is related with non-negative integer solutions of $n=\sum_{i\geq0}2^{i}N_{i}$.
If one expands some terms,
\begin{align*}
\,&\prod_{i\geq0}(\sum_{n\geq 0}p^{A}_{1}(n)q^{2^{i}n})\\
&=1+p^{A}_{1}(1)q+[p^{A}_{1}(2)+p^{A}_{1}(1)]q^2\\
&\;+[p^{A}_{1}(3)+p^{A}_{1}(1)p^{A}_{1}(1)]q^3\\
&\;+[p^{A}_{1}(4)+p^{A}_{1}(2)p^{A}_{1}(1)+p^{A}_{1}(2)+p^{A}_{1}(1)]q^4\\
&\;+[p^{A}_{1}(5)+p^{A}_{1}(1)p^{A}_{1}(2)+p^{A}_{1}(1)p^{A}_{1}(1)+p^{A}_{1}(3)p^{A}_{1}(1)]q^5\\&\quad\;\cdots
\end{align*}

\end{proof}

Now, we prove more general identities between $p^{A}_{}$ and $p^{A}_{\alpha}$. 

\begin{definition}
Let $n \in \mathbf{N}$ and $\{(a^{n,\alpha}_{11},a^{n,\alpha}_{12},a^{n,\alpha}_{13}\cdots)$, $(a^{n,\alpha}_{21},a^{n,\alpha}_{22},a^{n,\alpha}_{23}\cdots), \cdots \}$ be the set of all non-negative integer solutions of the indeterminate equation $n=N_{0}+(\alpha +1)N_{1}+(\alpha +1)^2N_{2}+\cdots=\sum_{i\geq 0}(\alpha +1)^{i}N_{i}$. In other words, 
\begin{align*}
n&=a^{n,\alpha}_{11}+(\alpha +1)a^{n,\alpha}_{12}+(\alpha +1)^2a^{n,\alpha}_{13}+\cdots\\
 &=a^{n,\alpha}_{21}+(\alpha +1)a^{n,\alpha}_{22}+(\alpha +1)^2a^{n,\alpha}_{23}+\cdots\\
 &\cdots
\end{align*}where $a^{n,\alpha}_{ij}\in\mathbf{N}$. Then the solution matrix of this equation is \[A_{n,\alpha}:=(a^{n,\alpha}_{ij}).\]
\end{definition}

\begin{theorem}
Let $\psi$ be an one to one function such that $\psi:\mathbf{N}\rightarrow \mathbf{N}$. And let $n\in\mathbf{N}$ and $A_{n,\alpha }=(a_{ij}^{n,\alpha})$ be the solution matrix of $n=\sum_{i\geq 0}(\alpha +1)^{i}N_{i}$. Then
\[
p^{A}_{}(n)=\sum_{i\geq 1}\prod_{j\geq 1}p^{A}_{\alpha}(a^{n,\alpha}_{ij})
\]
for all $\alpha\in\mathbf{N}\setminus\{0\}$.
\end{theorem}

\begin{proof}
It is well known fact that for $|q|<1$,

\begin{align*}
\sum_{n\geq 0}p^{A}_{\alpha}(n)q^{n}&=\prod_{n\geq 1}\frac{1-q^{(\alpha+1)\psi(n)}}{1-q^{\psi(n)}}\\
&=\prod_{n\geq 1}(1+q^{\psi(n)}+q^{2\psi(n)}+\cdots+q^{\alpha\psi(n)})
\end{align*}(see \cite{2}). On the other hand,

\begin{align*}\prod_{n\geq1}(1-q^{(\alpha+1)\psi(n)})&=\prod_{n\geq1}(1-q^{\psi(n)})\prod_{n\geq1}(1+q^{\psi(n)}+q^{2\psi(n)}+\cdots+q^{\alpha \psi(n)})\\
\prod_{n\geq1}(1-q^{(\alpha+1)^2\psi(n)})&=\prod_{n\geq1}(1-q^{\psi(n)})\prod_{n\geq1}(1+q^{\psi(n)}+q^{2\psi(n)}+\cdots+q^{\alpha \psi(n)})\\
&\,\times\prod_{n\geq1}(1+q^{(\alpha+1)\psi(n)}+q^{(\alpha+1)2\psi(n)}+\cdots+q^{(\alpha+1)\alpha\psi(n)})\\
\vdots\\
\prod_{n\geq1}(1-q^{(\alpha+1)^{I+1}\psi(n)})&=\prod_{n\geq1}(1-q^{\psi(n)})\\&\;\times\prod_{i=0}^{I}\prod_{n\geq1}(1+q^{(\alpha+1)^{i}\psi(n)}+\cdots+q^{(\alpha+1)^{i}\alpha \psi(n)})
\end{align*}
for all $I\in\mathbf{N}$. So,
\begin{align*}
1&=\lim_{I\rightarrow\infty}\prod_{n\geq1}(1-q^{(\alpha+1)^{I+1}\psi(n)})\\
&=\prod_{n\geq1}(1-q^{\psi(n)})\prod_{i=0}^{\infty}\prod_{n\geq1}(1+q^{(\alpha+1)^{i}\psi(n)}+\cdots+q^{(\alpha+1)^{i}\alpha \psi(n)}).
\end{align*}
Therefore,
\[
\prod_{n\geq1}\frac{1}{1-q^{\psi(n)}}=\prod_{i\geq0}\prod_{n\geq 1}(1+q^{(\alpha+1)^{i}\psi(n)}+\cdots+q^{(\alpha+1)^{i}\alpha\psi(n)})
\]
and
\[
\sum_{n\geq 0}p^{A}_{}(n)q^{n}=\prod_{i\geq 0}(\sum_{n\geq 0}p^{A}_{\alpha}(n)q^{(\alpha+1)^{i}n}).
\]
If one expands the above infinite product to an infinite series, one can express $p^{A}_{}(n)$ by $p^{A}_{\alpha}(n)$, $p^{A}_{\alpha}(n-1)$, $\cdots$, $p^{A}_{\alpha}(1)$ and one can see that the coefficient of $q^n$ is related with non-negative integer solutions of $n=\sum_{i\geq 0}(\alpha +1)^{i}N_{i}$ 
\end{proof}

\section{Inverse identities and some other similar identities}
 In section 1, we found the identities which express $p^{A}_{}(n)$ by $p^{A}_{\alpha}(n)$, $p^{A}_{\alpha}(n-1)$, $\cdots$, $p^{A}_{\alpha}(1)$. Now, we will find the inverse identities.
 
\begin{definition}
Let $E^{A}(n)$ be the number of  even partitions of $n$ without a restriction of the number of equal parts and $O^{A}(n)$ be the number of odd partitions of $n$ without a restriction of the number of equal parts. Then we define \[\bar{p}^{A}(n):=E^{A}(n)-O^{A}(n)\]
and $\bar{p}^{A}(0):=1.$
\end{definition}

\begin{definition}
Let $n=(\alpha+1)\sum_{i\geq0}2^{i}N_{i}$ be the indeterminate equation for $n, \alpha\in\mathbf{N}\setminus\{0\}$.
If this equation has solutions $(b^{n,\alpha}_{11}, b^{n,\alpha}_{12}, b^{n,\alpha}_{13}, \cdots), (b^{n,\alpha}_{21}, b^{n,\alpha}_{22}, b^{n,\alpha}_{23}, \cdots), \cdots$ where $b^{n,\alpha}_{ij}\in\mathbf{N}$, then we define the solution matrix of this equation by \[B_{n, \alpha}:=(b^{n,\alpha}_{ij}).\]
\end{definition}

\begin{definition}
Let $n, \alpha\in\mathbf{N}\setminus\{0\}$. If $n=(\alpha+1)\sum_{i\geq 0}2^{i}N_{i}$ has non-negative integer solutions, we define \[\Gamma^{\psi}_{\alpha}(n):=\sum_{i\geq 1}\prod_{j\geq 1}\bar{p}^{A}(b^{n,\alpha}_{ij})\] and if $n=(\alpha+1)\sum_{i\geq 0}2^{i}N_{i}$ does not have a non-negative integer solution, we define $\Gamma^{\psi}_{\alpha}(n):=0$. For $n=0$, we define $\Gamma^{\psi}_{\alpha}(0):=1$.
\end{definition}

\begin{theorem}
Let $\psi$ be an one to one function such that $\psi:\mathbf{N}\rightarrow \mathbf{N}$. Then 

\[
p^{A}_{\alpha}(n)=\sum^{n}_{i=0}p^{A}_{}(n-i) \Gamma^{\psi}_{\alpha}(i)
\]
for all $n, \alpha \in \mathbf{N}$.
\end{theorem}

\begin{proof}
For $|q|<1$,
\begin{align*}
\prod_{n\geq1}(1-q^{2(\alpha+1)\psi(n)})&=\prod_{n\geq 1}(1-q^{(\alpha+1)\psi(n)})\prod_{n\geq 1}(1+q^{(\alpha+1)\psi(n)})\\
\prod_{n\geq1}(1-q^{4(\alpha+1)\psi(n)})&=\prod_{n\geq 1}(1-q^{(\alpha+1)\psi(n)})\\
&\,\times\prod_{n\geq 1}(1+q^{(\alpha+1)\psi(n)})\prod_{n\geq 1}(1+q^{2(\alpha+1)\psi(n)})\\
\vdots\\
\prod_{n\geq1}(1-q^{2^{I+1}(\alpha+1)\psi(n)})&=\prod_{n\geq 1}(1-q^{(\alpha+1)\psi(n)})\prod_{i=0}^{I}\prod_{n\geq 1}(1+q^{2^i(\alpha+1)\psi(n)})
\end{align*}
for all $I\in\mathbf{N}$. So,
\begin{align*}
1&=\lim_{I\rightarrow\infty}\prod_{n\geq1}(1-q^{2^{I+1}(\alpha+1)\psi(n)})\\
&=\prod_{n\geq 1}(1-q^{(\alpha+1)\psi(n)})\prod_{i=0}^{\infty}\prod_{n\geq 1}(1+q^{2^i(\alpha+1)\psi(n)})
\end{align*}
and
\[
\prod_{n\geq 1}\frac{1-q^{(\alpha+1)\psi(n)}}{1-q^{\psi(n)}}
=\prod_{n\geq 1}\frac{1}{1-q^{\psi(n)}}\prod_{i\geq0}\prod_{n\geq 1}\frac{1}{1+q^{2^{i}(\alpha+1)\psi(n)}}.
\]On the other hand,
\[
\prod_{n\geq 1}\frac{1}{1+q^{2^{i}(\alpha+1)\psi(n)}}=\sum_{n\geq 0}\bar{p}^{A}(n)q^{2^{i}(\alpha+1)n}.
\]Therefore, if we define 

\[
\prod_{i\geq 0}\prod_{n\geq 1}\frac{1}{1+q^{2^{i}(\alpha+1)\psi(n)}}:=\sum_{n\geq 0}\Gamma^{\psi}_{\alpha}(n)q^{n},
\]then
\[
\Gamma^{\psi}_{\alpha}(n)=\sum_{i\geq 1}\prod_{j\geq 1}\bar{p}^{A}(b^{n,\alpha}_{ij})
\] when $n=(\alpha+1)\sum_{i\geq 0}2^{i}N_{i}$ has non-negative integer solutions and $\Gamma^{\psi}_{\alpha}(n)=0$ when $n=(\alpha+1)\sum_{i\geq 0}2^{i}N_{i}$ does not have a non-negative integer solution. \\Finally,
\begin{align*}
\sum_{n\geq 0}p^{A}_{\alpha}(n)q^{n}&=\prod_{n\geq 1}\frac{1}{1-q^{\psi(n)}}\prod_{i\geq 0}\prod_{n\geq 1}\frac{1}{1+q^{2^{i}(\alpha+1)\psi(n)}}\\
&=(\sum_{n\geq 0}p^{A}_{}(n)q^{n})(\sum_{n\geq 0}\Gamma^{\psi}_{\alpha}(n)q^{n})
\end{align*}and
\[
p^{A}_{\alpha}(n)=\sum^{n}_{i=0}p^{A}_{}(n-i) \Gamma^{\psi}_{\alpha}(i).
\]
\end{proof}

Next, we prove two similar theorems.
\begin{definition}
Let $E^{A}_{\alpha}(n)$ be the number of  even partitions of $n$ such that the number of equal parts is less than or equals to $\alpha\in \mathbf{N}$ and $O^{A}_{\alpha}(n)$ be the number of odd partitions of $n$ such that the number of equal parts is less than or equals to $\alpha\in\mathbf{N}\setminus\{0\}$. Then we define \[\bar{p}^{A}_{\alpha}(n):=E^{A}_{\alpha}(n)-O^{A}_{\alpha}(n)\]
and $\bar{p}^{A}_{\alpha}(0):=1.$
\end{definition}

\begin{theorem}
Let $A_n=(a^{n}_{ij})$ be the solution matrix of $n=\sum_{i\geq 0}2^{i}N_{i}$ and let $\psi$ be an one to one function such that $\psi:\mathbf{N}\rightarrow \mathbf{N}$. Then
\[
\bar{p}^{A}_{1}(n)=\sum_{i\geq 1}\prod_{j\geq 1}\bar{p}^{A}(a^{n}_{ij})
\]
for all $n\in\mathbf{N}$.
\end{theorem}

\begin{proof}
We proved that for $|q|<1,$
\[
1=\prod_{n\geq1}(1-q^{\psi(n)})\prod_{i=0}^{\infty}\prod_{n\geq1}(1+q^{2^{i}\psi(n)})
\]
in the proof of theorem 2.1.
So,
\[
\prod_{n\geq1}(1-q^{\psi(n)})=\prod_{i=0}^{\infty}\prod_{n\geq1}\frac{1}{(1+q^{2^{i}\psi(n)})}
\]
and
\[\sum_{n\geq0}\bar{p}^{A}_{1}(n)q^n=\prod_{i\geq0}(\sum_{n\geq0}\bar{p}^{A}(n)q^{2^{i}n}).\] This proves the theorem.
\end{proof}

\begin{theorem}
Let $A_{n, \alpha}=(a^{n, \alpha}_{ij})$ be the solution matrix of $n=\sum_{i\geq 0}(\alpha+1)^{i}N_{i}$ and let $\psi$ be an one to one function such that $\psi:\mathbf{N}\rightarrow \mathbf{N}$.
Then
\[
\bar{p}^{A}(n)=\sum_{i\geq 1}\prod_{j\geq 1}\bar{p}^{A}_{\alpha}(a^{n, \alpha}_{ij})
\]
for all $n\in\mathbf{N}$ and for all even natural number $\alpha$.
\end{theorem}

\begin{proof}
Let $\alpha$ be an even natural number. If $|q|<1,$
\begin{align*}
&\prod_{n\geq1}(1+q^{(\alpha+1)^{I+1}\psi(n)})\\
&=\prod_{n\geq1}(1+q^{\psi(n)})\prod_{i=0}^{I}\prod_{n\geq1}(1-q^{(\alpha+1)^{i}\psi(n)}+q^{(\alpha+1)^{i}2\psi(n)}-\cdots+q^{(\alpha+1)^{i}\alpha\psi(n)})
\end{align*}
for all $I\in\mathbf{N}$. So,
\begin{align*}
1&=\lim_{I\rightarrow\infty}\prod_{n\geq1}(1+q^{(\alpha+1)^{I+1}\psi(n)})\\
 &=\prod_{n\geq1}(1+q^{\psi(n)})\prod_{i=0}^{\infty}\prod_{n\geq1}(1-q^{(\alpha+1)^{i}\psi(n)}+q^{(\alpha+1)^{i}2\psi(n)}-\cdots+q^{(\alpha+1)^{i}\alpha\psi(n)}).
\end{align*}
Therefore,
\begin{align*}
\prod_{n\geq1}\frac{1}{(1+q^{\psi(n)})}&=\prod_{i=0}^{\infty}\prod_{n\geq1}(1-q^{(\alpha+1)^{i}\psi(n)}+q^{(\alpha+1)^{i}2\psi(n)}-\cdots+q^{(\alpha+1)^{i}\alpha\psi(n)})
\end{align*}
and
\[
\sum_{n\geq0}\bar{p}^{A}(n)=\prod_{i\geq0}(\sum_{n\geq0}\bar{p}^{A}_{\alpha}(n)q^{(\alpha+1)^{i}n}).
\]
This proves the theorem.
\end{proof}

\pagebreak

\section{Appendix: Examples}

In this appendix, we will consider two identities for $p^{A}_{}(10)$ and $p^{A}_{1}(10)$.

I) Let us consider the expression of $p^{A}_{}(10)$ by $p^{A}_{1}$s. The indeterminate equation for this identity is
\[
10=N_0+2N_1+4N_2+\cdots.
\]The solution matrix of this equation and the identity for $p^{A}_{}(10)$ are

$\left(
 \begin{matrix}
 0 & 1 & 0 & 1 & \cdots \\ 
 0 & 1 & 2 & 0 & \cdots \\ 
 0 & 3 & 1 & 0 & \cdots \\ 
 0 & 5 & 0 & 0 & \cdots \\ 
 2 & 2 & 1 & 0 & \cdots \\ 
 2 & 4 & 0 & 0 & \cdots \\ 
 2 & 0 & 0 & 1 & \cdots \\ 
 2 & 0 & 2 & 0 & \cdots \\ 
 4 & 1 & 1 & 0 & \cdots \\ 
 4 & 3 & 0 & 0 & \cdots \\ 
 6 & 2 & 0 & 0 & \cdots \\ 
 6 & 0 & 1 & 0 & \cdots \\ 
 8 & 1 & 0 & 0 & \cdots \\ 
 10 & 0 & 0 & 0 & \cdots
 \end{matrix}\right) 
 $,$\quad$ $\begin{array}{l}
 p^{A}_{}(10) \\
  \\
 =p^{A}_{1}(1)p^{A}_{1}(1)+p^{A}_{1}(1)p^{A}_{1}(2)+p^{A}_{1}(1)p^{A}_{1}(3) \\
  +p^{A}_{1}(5)+p^{A}_{1}(2)p^{A}_{1}(2)p^{A}_{1}(1)\\ 
  +p^{A}_{1}(2)p^{A}_{1}(4)+p^{A}_{1}(2)p^{A}_{1}(1) \\ 
 +p^{A}_{1}(2)p^{A}_{1}(2)+p^{A}_{1}(4)p^{A}_{1}(1)p^{A}_{1}(1) \\ 
 +p^{A}_{1}(4)p^{A}_{1}(3)+p^{A}_{1}(6)p^{A}_{1}(2)+p^{A}_{1}(6)p^{A}_{1}(1) \\ 
 +p^{A}_{1}(8)p^{A}_{1}(1)+p^{A}_{1}(10).
 \end{array} 
 $

 Now, we calculate $p^{A}_{}(10)$ for three set of parts.\\1) $\psi(\mathbf{N})=\{p\;|\;p\;is\;a\;prime\}$\\Let us calculate partition numbers.\\
 $
 \begin{array}{ll}
 10= & 5+5 \\ 
  & 3+2+5 \\ 
  & 3+3+2+2 \\ 
  & 7+3 \\ 
  & 2+2+2+2+2
 \end{array} 
 $, $
 \begin{array}{ll}
 p^{A}_{1}(1)=0 & p^{A}_{1}(5)=2 \\ 
 p^{A}_{1}(2)=1 & p^{A}_{1}(6)=0 \\ 
 p^{A}_{1}(3)=1 & p^{A}_{1}(8)=1 \\ 
 p^{A}_{1}(4)=0 & p^{A}_{1}(10)=2.
 \end{array} 
 $\\So, $p^{A}_{}(10)=5$ and since $p^{A}_{1}(1)=p^{A}_{1}(4)=p^{A}_{1}(6)=0$, 
 \[
 p^{A}_{}(10)=p^{A}_{1}(5)+p^{A}_{1}(2)p^{A}_{1}(2)+p^{A}_{1}(10)=2+1+2=5.
 \]
\\2) $\psi(\mathbf{N})=\{n^2\;|\;n\in\mathbf{N} \} $.\\Let us calculate the partition numbers.\\
$
 \begin{array}{ll}
 10= & 1+9 \\ 
  & 1+1+1+1+1+1+1+1+1+1 \\ 
  & 4+1+1+1+1+1+1 \\ 
  & 4+4+1+1
 \end{array} 
$, $
 \begin{array}{ll}
 p^{A}_{1}(1)=1 & p^{A}_{1}(5)=1 \\ 
 p^{A}_{1}(2)=0 & p^{A}_{1}(6)=0 \\ 
 p^{A}_{1}(3)=0 & p^{A}_{1}(8)=0 \\ 
 p^{A}_{1}(4)=1 & p^{A}_{1}(10)=1.
 \end{array} 
$\\So, $p^{A}_{}(10)=4$ and since $p^{A}_{1}(2)=p^{A}_{1}(3)=p^{A}_{1}(6)=p^{A}_{1}(8)=0$,
\begin{align*}
p^{A}_{}(10)&=p^{A}_{1}(1)p^{A}_{1}(1)+p^{A}_{1}(5)+p^{A}_{1}(4)p^{A}_{1}(1)p^{A}_{1}(1)+p^{A}_{1}(10)\\
&=1+1+1+1=4.
\end{align*}

3) $\psi(\mathbf{N})=\{n\;|\;n\;is\;an\;odd\;number\}$.\\Let us calculate the partition numbers.\\
$
 \begin{array}{ll}
 10= & 1+1+1+1+1+1+1+1+1+1 \\ 
  & 3+1+1+1+1+1+1+1 \\ 
  & 5+1+1+1+1+1 \\
  & 7+1+1+1 \\
  & 9+1 \\
  & 3+3+1+1+1+1 \\
  & 3+3+3+1 \\
  & 5+5 \\
  & 7+3 \\
  & 3+5+1+1
 \end{array} 
$, $
 \begin{array}{ll}
 p^{A}_{1}(1)=1 & p^{A}_{1}(5)=1 \\ 
 p^{A}_{1}(2)=0 & p^{A}_{1}(6)=1 \\ 
 p^{A}_{1}(3)=1 & p^{A}_{1}(8)=2 \\ 
 p^{A}_{1}(4)=1 & p^{A}_{1}(10)=2.
 \end{array}$\\So, $p^{A}_{}(10)=10$ and since $p^{A}_{1}(2)=0$,\\
 
 \begin{align*}
  p^{A}_{}(10) & =p^{A}_{1}(1)p^{A}_{1}(1)+p^{A}_{1}(1)p^{A}_{1}(3)+p^{A}_{1}(5) \\
  & \,\,+p^{A}_{1}(4)p^{A}_{1}(1)p^{A}_{1}(1)+p^{A}_{1}(4)p^{A}_{1}(3) \\
  & \,\,+p^{A}_{1}(6)p^{A}_{1}(1)+p^{A}_{1}(8)p^{A}_{1}(1)+p^{A}_{1}(10) \\
  & \,\,=1+1+1+1+1+1+2+2 \\
  & \,\,=10. 
 \end{align*}
 
 II) Let us consider the identity for $p^{A}_{1}(10)$. The indeterminate equation for this identity is \[n=2N_0+4N_1+8N_2+\cdots\] and solutions of this equation are
\[
\begin{tabular}{|ll|}
\hline n & solutions \\ 
\hline 2 & $(1, 0, 0, \cdots)$ \\ 
\hline 4 & $(2, 0, 0, \cdots)$, $(0, 1, 0, \cdots)$ \\ 
\hline 6 & $(3, 0, 0, \cdots)$, $(1, 1, 0, \cdots)$ \\ 
\hline 8 & $(4, 0, 0, \cdots)$, $(2, 1, 0, \cdots)$, $(0, 0, 1, \cdots)$, $(0, 2, 0, \cdots)$ \\ 
\hline 10 & $(5, 0, 0, \cdots)$, $(1, 2, 0, \cdots)$, $(1, 0, 1, \cdots)$, $(3, 1, 0, \cdots)$ \\ 
\hline 
\end{tabular}
\]and $\emptyset$ for $n=odd\;number$.
So,
\[
\begin{array}{ll}
n & \Gamma^{\psi}_{1}(n) \\ 
0 & 1 \\ 
2 & \bar{p}^{A}(1) \\ 
4 & \bar{p}^{A}(2)+\bar{p}^{A}(1) \\ 
6 & \bar{p}^{A}(3)+\bar{p}^{A}(1)\bar{p}^{A}(1) \\ 
8 & \bar{p}^{A}(4)+\bar{p}^{A}(2)\bar{p}^{A}(1)+\bar{p}^{A}(2)+\bar{p}^{A}(1) \\ 
10 & \bar{p}^{A}(5)+\bar{p}^{A}(2)\bar{p}^{A}(1)+\bar{p}^{A}(1)\bar{p}^{A}(1)+\bar{p}^{A}(3)\bar{p}^{A}(1) 
\end{array}\]
and $\Gamma^{\psi}_{1}(n)=0$ for $n=odd\,number$.\\Therefore, if we denote $p^{A}_{}(n)$ by $p^{A}(n)$ and $\bar{p}^{A}(n)$ by $\bar{p}^{A}(n)$, 
\[\,\]
$
\begin{array}{ll}
p^{A}_{1}(10) & =p^{A}(10)+p^{A}(8)\bar{p}^{A}(1)+p^{A}(6)[\bar{p}^{A}(2)+\bar{p}^{A}(1)] \\
 & \,\,+p^{A}(4)[\bar{p}^{A}(3)+\bar{p}^{A}(1)\bar{p}^{A}(1)] \\  & \,\,+p^{A}(2)[\bar{p}^{A}(4)+\bar{p}^{A}(2)\bar{p}^{A}(1)+\bar{p}^{A}(2)+\bar{p}^{A}(1)] \\
 & \,\,+\bar{p}^{A}(5)+\bar{p}^{A}(2)\bar{p}^{A}(1)+\bar{p}^{A}(1)\bar{p}^{A}(1)+\bar{p}^{A}(3)\bar{p}^{A}(1).
\end{array}
$
\[\,\]
Now, we claulate for two set of parts $\{p\;|\;p\;is\;a\;prime\}$ and $\{n^{2}\;|\;n\in\mathbf{N}\}$.\\
1) For $\{p\;|\;p\;is\;a\;prime\}$, $p^{A}_{1}(10)=2$ and\\
$
\begin{array}{ll}
p^{A}(2)=1 & \bar{p}^{A}(1)=0 \\ 
p^{A}(4)=1 & \bar{p}^{A}(2)=-1 \\ 
p^{A}(6)=2 & \bar{p}^{A}(3)=-1 \\ 
p^{A}(8)=3 & \bar{p}^{A}(4)=1 \\ 
p^{A}(10)=5 & \bar{p}^{A}(5)=0. 
\end{array}$
\\If we calculate $p^{A}_{1}(10)$,
\begin{align*}
 p^{A}_{1}(10) & =2 \\ 
 & =5+3\times0+2\times[(-1)+0]+1\times[(-1)+0\times0] \\ 
 & \,\,+1\times[1+(-1)\times0+(-1)+0]+0\\
 &\,\,+0\times(-1)+0\times0+0\times(-1).
\end{align*}

2) For $\{n^{2}\;|\;n\in\mathbf{N}\}$, $p^{A}_{1}(10)=1$ and\\
$
\begin{array}{ll}
p^{A}(2)=1 & \bar{p}^{A}(1)=-1 \\ 
p^{A}(4)=2 & \bar{p}^{A}(2)=1 \\ 
p^{A}(6)=2 & \bar{p}^{A}(3)=-1 \\ 
p^{A}(8)=3 & \bar{p}^{A}(4)=0 \\ 
p^{A}(10)=4 & \bar{p}^{A}(5)=0. 
\end{array}$
\\If we calculate $p^{A}_{1}(10)$,
\begin{align*}
 p^{A}_{1}(10) & =1 \\
 & =4+3\times(-1)+2\times[1+(-1)] \\ 
 & \,\,+2\times[(-1)+(-1)\times(-1)] \\ 
 & \,\,+1\times[0+1\times(-1)+1+(-1)] \\
 & \,\,+0+(-1)\times1+(-1)\times(-1)+(-1)\times(-1).
\end{align*}


\begin{thebibliography}{widestlabel}
\bibitem{1} Hua, L.-K. {\it Introduction to number theory}, (Trans. from Chinese by Peter Shiu), Springer-Verlag, Berlin, 1982.
\bibitem{2} Rose, H. E. {\it A course in number theory}, 2nd edition, Oxford Science Publications, Oxford, 1994.
\end{thebibliography}
\end{document}